\documentclass{amsart}

\newtheorem{theorem}{Theorem}[section]
\newtheorem{lemma}[theorem]{Lemma}

\newtheorem{proposition}[theorem]{Proposition}

\theoremstyle{definition}
\newtheorem{definition}[theorem]{Definition}

\theoremstyle{remark}
\newtheorem{remark}[theorem]{Remark}

\usepackage{xy}
\xyoption{all}

\usepackage{amsmath}
\usepackage{amsfonts}
\usepackage{amssymb}

\def\sideremark#1{\ifvmode\leavevmode\fi\vadjust{\vbox to0pt{\vss % the remark
      \hbox to 0pt{\hskip\hsize\hskip1em           %                will appear only
 \vbox{\hsize3cm\tiny\raggedright\pretolerance10000%                on the side
 \noindent #1\hfill}\hss}\vbox to8pt{\vfil}\vss}}}%
\newcommand{\bbR}{\mathbb{R}}
\newcommand{\bbC}{\mathbb{C}}
\newcommand{\bbZ}{\mathbb{Z}}
\newcommand{\sog}{\mathbf{SO}}

\newcommand{\soa}{\mathfrak{so}}

\newcommand{\sph}{\mathbb{S}}

\newcommand{\st}{V}
\newcommand{\stb}{\mathcal{V}}
\newcommand{\br}{B}
\newcommand{\brb}{\mathcal{B}}
\newcommand{\pb}{\mathcal{P}}
\newcommand{\obs}{\mathfrak{o}}
\newcommand{\bbP}{\mathbb{P}}
\newcommand{\scm}{{\scriptscriptstyle -}}
\newcommand{\ten}{\Upsilon}

\begin{document}
\title[The topological obstructions to the irreducible SO(3) structure]{The topological obstructions to the existence of an irreducible SO(3) structure on a five manifold}
%\primaryclass{57R15;53C10}
%53B15 (local DG, other connections)
%53B50 (loc DG applic to physics)
%53C10 (glob DG, G-str)
%53C15 (glob DG, general geom str on mnfolds)
%57R15 (diff top, specialized structures on manifoilds(spin))
%53C30 (homogeneous manifolds)
%\vskip 1.truecm 
\author{Marcin Bobie\'nski} \address{Instytut Matematyki,
Universytet Warszawski, ul. Banacha 2, Warszawa, Poland}
\email{mbobi@mimuw.edu.pl} \thanks{This research was supported by
the KBN Grant No 1 P03A 015 29.\\
During the preparation of this
article the author was member of the VW Junior Research Group ``Special
Geometries in Mathematical Physiscs'' at Humboldt University in
Berlin.}
\date{\today}

\begin{abstract}
A nonstandard (maximal) inclusion $\sog(3)\subset \sog(5)$ associated with the irreducible representation $\rho_5$ of $\sog(3)$ in $\bbR^5$ is considered. The topological obstructions for admitting the $\sog(3)$ structure on the frame bundle over 5-manifold are investigated. The necessary and sufficient conditions are formulated.
\end{abstract}
\maketitle
\section{Introduction}

In Ref. \cite{bn} we introduced and investigated the irreducible (maximal) $\sog(3)$ structure on a 5-dimensional manifold $M$. We found and described the tensor object reducing the structure group of the frame bundle of $M$ from $\sog(5)$ to the irreducible $\sog(3)$. That paper was mainly devoted to the local analysis of the geometry of manifold $M$ equipped with such a structure. Our motivation for investigation of structures of such kind was the paper \cite{Fried} of Th. Friedrich, where he listed especially interesting types of special geometries in low dimensions. There are such interesting geometries like $\mathbf{G_2}$ structure in dimension 7, $\mathbf{Spin}(7)$ structure in dimension 8, $\mathbf{Spin}(9)$ structure in dimension 16, $\mathbf{F_4}$ structure in dimension 26; Friedrich also adds the $\sog(3)$ structure in dimension 5 to this list.

In the case of any structure it is interesting to know under which topological conditions the structure exists on a manifold $M$. For example, it is well known that the $\mathbf{Spin}$ structure on an oriented manifold $M$ does exist provided the second Stiefel-Whitney class of the tangent bundle $w_2(TM)$ vanishes. 

The main goal of this paper is to prove the following criterion. There exists the maximal $\sog(3)$ structure on an oriented 5-dimensional manifold $M$ if and only if there exists the standard $\sog(3)$ structure (i.e.\ the tangent bundle decomposes $TM=E^3\oplus\theta^2$) and the first Pontryagin class $p_1(TM)$ is divisible by 5. This result is used to construct non-trivial examples of 5-manifolds equipped with the maximal $\sog(3)$ structure -- see Proposition \ref{pr:cs} below.

\subsection{The irreducible $\sog(3)$ structure}
\label{sec:so3str}

To fix the notation let us recall the explicit construction of the unique 5-dimensional representation of $\sog(3)$. We identify $\bbR^5$ with the subspace of $3\times 3$ real matrices
\begin{equation}
\mathbb{M}^5=\{ A\in M_{3\times 3}(\mathbb{R}):~~ A^{\rm T}=A,~~{\rm tr}(A)=0~\}. \label{r5}
\end{equation}
The action of $\sog(3)$ on $\mathbb{M}^5$ is given by
\begin{equation}
\rho_5(h) A~=~h~A~h^{\rm T},\quad\quad\quad\forall~h\in\sog(3),~~~~~~~A\in \mathbb{M}. \label{so35}
\end{equation}
Each $\rho_5(h)$ defines the orthogonal transformation of $\mathbb{M}^5$, so the representation $\rho_5$ defines the inclusion
\begin{equation}
\iota_5: \sog(3)\hookrightarrow\sog(5), \label{iota}
\end{equation}
which is essentially different from the standard (diagonal) inclusion 
\begin{equation}
  \label{ist}
  j:  \sog(3)\hookrightarrow\sog(5).
\end{equation}
The image $\iota_5(\sog(3))$ will be called \emph{irreducible} or \emph{maximal} $\sog(3)$ and the usual one $j(\sog(3))$ will be called \emph{standard} $\sog(3)$.

The irreducible $\sog(3)$ structure on a 5-dimensional manifold $M$ is the reduction of the structure group of the frame bundle $FM$ to the irreducible $\sog(3)$. It was proved in \cite{bn} that this reduction is obtained by a pair of tensors $(g,\ten)$ satisfying the following relations.
\begin{definition}
\label{df:so35}
The \emph{irreducible} $\sog(3)$ structure on a 5-dimensional manifold $M$ is a triple $(M,g,\ten)$, where $g$ is a Riemannian metric and $\ten$ is a rank 3 tensor field defining vector bundle morphism
\[
\ten: TM \longrightarrow \mathrm{End}(TM),\qquad v\mapsto \ten_v\in \mathrm{End}(TM).
\]
This morphism satisfies the following conditions
\begin{enumerate}
\item it is totally symmetric, i.e.\ $g(u,\ten_v w) = g(w,\ten_v u) = g(u,\ten_w v)$ for all $u,v,w\in TM$;
  \item it is trace free $\mathrm{tr}(\ten_v)=0$,
  \item for any vector v $\in\bbR^5$ the following identity holds
\begin{equation}
\label{tid}
\ten_v^2 v = g(v,v) v.
\end{equation}
\end{enumerate}
\end{definition}

\begin{remark}
Let me recall briefly how we obtain the tensor $\ten$ from the irreducible $\sog(3)$ structure (see \cite{bn} for details). The explicit realization of the irreducible representation $\rho_5$ is constructed via identification of the point in $\bbR^5$ with the symmetric, trace-free $3\times 3$ matrix $A$. The group action is realized by the adjoint transformation. Thus, the determinant $\det A$ is homogeneous invariant polynomial of degree 3; it defines the rank 3 symmetric tensor $\ten$.
\end{remark}

\begin{remark}
It is worth to note that the tensor $\ten$ alone suffices to reduce the structure group to the irreducible $\sog(3)$. One can read out the Riemannian metric from the identity (\ref{tid}). The alternative definition of the irreducible $\sog(3)$ structure on a manifold $M$, involving only the bundle morphism $\ten$, is the following.
\begin{enumerate}
\item For any vector $v$ the endomorphism $\ten_v$ is trace-free.
\item Any vector $v$ is an eigenvector of $\ten_v^2$. The respective eigenvalue form $g(v,v)$ is quadratic and we assume to be positively defined.
\item The positive, quadratic form, defined in the previous point, provides a Riemannian metric $g$; the tensor $\ten$ is totally symmetric with respect to this metric.
\end{enumerate}
\end{remark}

Now, the problem of existence of the irreducible $\sog(3)$ structure arises. The main result of this paper is the following theorem.

\begin{theorem}
\label{th:top}
Let $M$ be an orientable 5-dimensional manifold. There exists an irreducible $\sog(3)$ structure on $M$ if and only if
\[
TM = E^3\oplus \theta^2,\qquad p_1(TM)= 5\; \tilde{p},\quad \tilde{p}\in H^4(M;\bbZ),
\]
where $\theta^2$ is the 2-dimensional trivial bundle and $p_1(TM)$ is the Pontrjagin class.
\end{theorem}

The above theorem provides wide variety of non-trivial compact examples of manifolds admitting the irreducible $\sog(3)$ structure. 

\begin{proposition}
\label{pr:cs}
Let $S$ be a complex surface and $M=S\times \sph^1$, where $\sph^1$ is a circle.
\begin{enumerate}
\item \label{it:st} There exist the standard $\sog(3)$ structure on $M$ if and only if 
\[
\chi(S) \equiv 0 \mod 2.
\]
\item \label{it:br} There exist the irreducible $\sog(3)$ structure on $M$ if and only if 
\[
\chi(S)\equiv 0\mod 2, \quad \text{and} \quad \sigma(S) \equiv 0\mod 5,
\]
where $\sigma(S)$ is the signature of $S$.
\end{enumerate}
\end{proposition}
\begin{proof}~

\emph{Point (\ref{it:st}).} We fix an arbitrary Riemannian metric on the tangent bundle $TM$ and we consider the corresponding Stiefel fibration $\stb\rightarrow M$. Its fibre $\st$ is the Stiefel manifold consisting of pairs $(v_1,v_2)$ of unit, orthogonal tangent vectors $v_1,v_2\in T M$. The standard $\sog(3)$ structure on $M$ corresponds to a section of $\stb$. The only obstruction in the construction of section over $S$ is the 4-th Stiefel-Whitney class \cite{milch}
\[
w_4(TM|_S) = w_4(TS).
\]
The evaluation of the latter class $w_4(TS)$ on the fundamental cycle of $S$ (i.e.\ ``integration'' over $S$) gives (see \cite{milch} again)
\[
<w_4(TS),S> \equiv \chi(S) \mod 2.
\]
Thus the Stiefel fibration restricted to $S$, $\stb|_S\rightarrow S$ admits a section $f$ if and only if $\chi(S)$ is even. The section $f\in\Gamma(\stb|_S\rightarrow S)$ can be further prolongated on the whole $M=S\times \sph^1$. This finishes the proof of point (\ref{it:st}).
%of Proposition \ref{pr:cs}.

The second point of the Proposition is the consequence of Theorem \ref{th:top} combined with the Hirzebruch signature formula $\sigma(S) = 3 <p_1(TS),S>$.\\
\end{proof}

\begin{remark}
The above proposition implies, for example, that there exist the irreducible $\sog(3)$ structure on $\widetilde{\bbC\bbP^2}\times\sph^1$, where $\widetilde{\bbC\bbP^2}$ is the projective space $\bbC\bbP^2$ with one point blown up.
\end{remark}

\section{Proof of Theorem \ref{th:top}}

%We denote the 7-dimensional homogeneous spaces
We consider the following 7-dimensional homogeneous spaces
\[
\st = \sog(5)/j(\sog(3)), \qquad \br = \sog(5)/\iota_5(\sog(3)),
\]
which are the quotients of $\sog(5)$ by the standard and the irreducible $\sog(3)$ respectively. The homogeneous space $\st=V_{2,5}$ is the Stiefel manifold of pairs of unit orthogonal vectors in $\bbR^5$. The second space $\br$ is known as the Berger space \cite{ks,gks}.

In the following proposition we recall one result about homologies of the homogeneous spaces $\st$ and $\br$.
\begin{proposition}[\cite{ks,hus}]
\label{pr:UVtop}
The homogeneous spaces $\st,\ \br$ are simply-connected manifolds with the following integral homologies
\begin{align*}
  H_3(\st) = \bbZ_2, \qquad H_3(\br) =\bbZ_{10},\\
  H_i (\st) = H_i(\br) = \bbZ, \qquad i=0,7,\\
%  H_3(\st) = \bbZ_2, \qquad H_3(\br) =\bbZ_{10},\\
  H_i (\st) = H_i(\br) = 0, \qquad i \neq 0,3,7.
\end{align*}
%and all other homology groups vanish.
\end{proposition}
%\medskip

\begin{proof}[Proof of Theorem \ref{th:top}]
The main idea of proof goes as follows. We consider the $\sog(5)$ principal bundle and the associated Berger and Stieffel fibrations. There does not exist the morphism of the latter fibrations. Nevertheless, we consider the Postnikov towers of fibrations and we construct a morphism of these towers up to the 4-th stages. Since the dimension of the base (i.e.\ the manifold $M$) is 5, this construction efficiently mimic the fibrations morphism.

%The main idea of the proof is to find a replacement for the $\sog(5)$-equivariant morphism $\br \rightarrow \st$ inducing the projection of the homology group $\bbZ_{10}=\bbZ_2\oplus\bbZ_5\rightarrow \bbZ_2$ (such a map obviously does not exist). The replacement is the morphism of the Postnikov towers up to the 4-th stage. 

Let us recall (see \cite{span}) the notation related to the Postnikov towers of fibration. A fibration $E\rightarrow M$ provides a tower of fibre spaces $p_j:E_j\rightarrow E_{j-1}$ with fibers being the Eilenberg-McLane spaces $K(A_j,n_j)$ and with the following commutative diagram
\[
\xymatrix{
& K(A_3,n_3) \ar[d] & K(A_2,n_2) \ar[d] & K(A_1,n_1) \ar[d] &  \\
\cdots  \ar[r] & E_3 \ar[r]^{p_3} & E_2 \ar[r]^{p_2} & E_1 \ar[r]^{p_1} & M\\
& & E \ar[ul] \ar[u] \ar[ur] \ar[rr]^{p} & &  M.
}
\]

Let $\sog(5)\rightarrow \pb \rightarrow M$ be an $\sog(5)$-principal bundle over a base $M$. Let 
\[
\brb = \pb\times_{\sog(5)} \br,\qquad \stb = \pb\times_{\sog(5)} \st
\]
be associated bundles with fibers $\br$ and $\st$ respectively. 

\begin{lemma}
\label{lem:post}
There exist fibration morphisms $\phi_3,\, \phi_4$, between the respective Postnikov fibrations of $\brb$ and $\stb$, which induce isomorphisms of fibers and make the following diagram commutative
\[
\xymatrix{
 & \vdots \ar[d] & \vdots \ar[d] & \\
K(\bbZ_2, 4) \ar[r] &\brb_4 \ar[r]^{\phi_4} \ar[d] & \stb_4 \ar[d]& \ar[l] K(\bbZ_2, 4)\\
K(\bbZ_2, 3) \ar[r] &\brb_3^{(2)} \ar[r]^{\phi_3} \ar[d] & \stb_3 \ar[d]& \ar[l] K(\bbZ_2, 3)\\
K(\bbZ_5, 3) \ar[r] &\brb_3^{(1)} \ar[r]^{p_1} \ar[d]^{p_1} & M\\
& M.&  & \\
}
\]
\end{lemma}

Now let us consider an oriented 5-dimensional manifold $M$. We have the fibrations $\stb$ and $\brb$ defined above.
%denote by $\stb$ the fibration of Stiefel manifolds associated to the tangent bundle; the fibration of Berger spaces is denoted by $\brb$. 
The existence of the usual $\sog(3)$ and the irreducible $\sog(3)$ structure is equivalent to the existence of a section of $\stb$ and $\brb$ respectively. 

The existence of a section of $\brb$ is equivalent to the existence of a section over $M$ of the 4-th stage in the Postnikov tower $\brb_4$. The latter statement is a consequence of fact that the base $M$ is 5-dimensional and the fibers of higher fibrations in the Postnikov tower are 4-connected, so the construction of section over $M$ of a higher fibration is unobstructed. The analogous statement is true for the Stieffel fibration $\stb$.

%The existence of these sections are equivalent to the existence of sections over $M$ of the 4-th stages of the Postnikov towers $\stb_4$ and $\brb_4$. The latter statement is a consequence of the fact that the base is 5-dimensional and the fibers of higher fibrations in the Postnikov towers are 4-connected, so the construction of section over $M$ of higher fibrations is unobstructed.

Using the Lemma \ref{lem:post}, we deduce that there is one to one correspondence between sections of $\brb_4$ and the following pairs of sections
\[
\Gamma(\brb_4 \to M) \stackrel{1:1}{\longleftrightarrow}  \Big\{\Gamma(\brb_3^{(1)} \to M),\quad\Gamma(\stb_4 \to M)  \Big\}.
\]
The obstruction to the section of $\brb_3^{(1)} \to M$ is the characteristic class $\obs\in H^4(M;\bbZ_5)$. Analyzing the cohomology groups of oriented grassmannian \cite{milch}, we deduce that it must be
\begin{equation}
\label{obst}
\obs = \lambda\; [p_1(TM)]_5,\quad \lambda \in\bbZ_5,
\end{equation}
where $[p_1(TM)]_5$ is the reduction modulo 5 of the Pontrjagin class. Thus, we have the alternative: either $\lambda\neq 0$ or any oriented 5-manifold $M$ such that $TM=E^3\oplus\theta^2$ admits the irreducible $\sog(3)$ structure. 

To finish the proof we show the following
\begin{lemma}
\label{lem:5pon}
Let $\pb$ be an $\sog(3)$ principal bundle and $\rho_5$ be the irreducible 5-dimensional representation of $\sog(3)$. Then, the first Pontrjagin class of the associated bundle is divisible by 5 i.e.\ 
\[
p_1(\pb\times_{\rho_5} \bbR^5) \equiv 0 \mod 5.
% 5\; \tilde{p},\quad \tilde{p}\in H^4(M;\bbZ),
\]
\end{lemma}

Let $K3$ denote the K3-surface. Let us recall that $\chi(K3)=24$ and $\sigma(K3)=-16$. The manifold $M=K3\times \sph^1$ provides an example of 5-manifold whose tangent bundle decomposes $TM=E^3\oplus\theta^2$ (see Proposition \ref{pr:cs}) and whose Pontrjagin class $<p_1(TM),K3>=-3\cdot 16$ is \emph{not} divisible by 5. 
This shows that $\lambda \neq 0$ in (\ref{obst}) and hence the Theorem \ref{th:top} is proved.\\
%Thus, $\lambda \neq 0$ and the Theorem \ref{th:top} is proved.\\
\end{proof}

\begin{remark}
Actually we have proven that the statement of the theorem remains valid after replacing the tangent bundle to a manifold with an arbitrary, oriented bundle of rank 5 over a base of dimension at most 5.
\end{remark}
%\medskip

\begin{proof}[Proof of Lemma \ref{lem:post}]
The fibers of $p_1$ are Eilenberg-McLane spaces $K(\bbZ_5,3)$ whose cohomologies with $\bbZ_2$ coefficients are trivial, at least up to gradation 5. Thus the map $p_1$ induces an isomorphism of cohomologies with $\bbZ_2$ coefficients (up to 5-th gradation). The fibrations $\brb_3^{(2)} \rightarrow \brb_3^{(1)}$ and $\stb_3\rightarrow M$ are determined by the characteristic class of the respective fibration. In both cases $\br$ and $\st$ the generator of the 3-rd homotopy group is the image of the generator of $\pi_3(\sog(5))=\bbZ$ under the canonical projections $\sog(5)\to \br$ and $\sog(5)\to \st$. Thus the characteristic classes in both cases of Berger and Stiefel fibrations coincide. It it known \cite{milch} that in the case of Stiefel fibration this obstruction class is equal to the 4-th Stiefel-Whitney class $w_4(TM)$. 

Summing up the above considerations, we have proved that the fibration $\brb_3^{(2)} \rightarrow \brb_3^{(1)}$ is isomorphic to the pull-back via $p_1$ of the $\stb_3$ fibration i.e.\ $p_1^*(\stb_3\rightarrow M)$. Thus, there exists the morphism $\phi_3$ as in the lemma.

To construct the morphism $\phi_4$ we show that the following fibrations over $\brb_3^{(2)}$ are isomorphic:
\[
(\brb_4 \rightarrow \brb_3^{(2)}) \cong \phi_3^* (\stb_4 \rightarrow \stb_3).
\]
These fibrations are determined by the second characteristic elements $k_\brb\in H^5(\brb_3^{(2)};\bbZ_2)$ and $k_\stb\in H^5(\st_3;\bbZ_2)$ respectively; so, it is enough to prove that $\phi_3^*(k_\stb) =k_\brb$. Since the bundles $\brb$ and $\stb$ are constructed out of the $\sog(5)$ principal bundle, it is enough to determine these characteristic elements for the tautological bundle over the classifying space $\mathbf{BSO(5)}$ i.e.\ the oriented grassmannian 
\[
M=G_5(\bbR^\infty).
\]
We consider the spectral sequences of the following fibrations
\[
\brb\rightarrow \brb_3^{(1)}, \quad \brb_3^{(2)}\rightarrow \brb_3^{(1)}, \quad \stb\rightarrow M,\quad \stb_3\rightarrow M.
\]
We omit the details of calculations which are quite standard. In this calculations we use the following known facts (with $\bbZ_2$ coefficients assumed) \cite{motan}.
\begin{enumerate}
\item The Steenrod algebra structure of cohomologies of the Eilenberg-McLane space $K(\bbZ_2,3)$.
\item The non trivial cohomologies of the Berger space $\br$ and the Stiefel manifold $\st$ are located in gradation 3 and 4. The Steenrod operation $Sq^1$ gives rise to the isomorphism of these spaces. The latter property is the consequence of the fact that $Sq^1$ coincides to the Bockstein homomorphism.
\item The transgression operation $\tau$ in the spectral sequence commutes with the Steenrod squaring operations $Sq^j$.
\item The characteristic elements $k_\brb$ and $k_\stb$ map to zero in cohomologies of the total spaces of fibrations, i.e.\ $H^5(\brb;\bbZ_2)$ and $H^5(\stb;\bbZ_2)$ respectively.
\end{enumerate}
Using the above facts we can unequally determine characteristic elements $k_\brb$ and $k_\stb$. Since the Steenrod algebra structure of cohomologies of the Berger space $\br$ and the Stiefel manifold $\st$ coincide, the respective characteristic elements also coincide; more precisely, they are related by $\phi_3^*$ which gives rise to the isomorphism of 5-th cohomologies (with $\bbZ_2$ coefficients)
\[
k_\brb =\phi_3^* k_\stb.
\]
Thus the morphism $\phi_4$ can be constructed.\\
\end{proof}

\begin{proof}[Proof of Lemma \ref{lem:5pon}]
Let $\pb$ be a principal $\sog(3)$ bundle, $\rho_3$ and $\rho_5$ be irreducible representation of $\sog(3)$ in dimension 3 and 5 respectively. We show the following relation among the Pontrjagin classes of the associated bundles
\begin{equation}
  \label{ponrel}
  p_1(\pb\times_{\rho_5} \bbR^5) = 5\cdot p_1(\pb\times_{\rho_3} \bbR^3).
\end{equation}

Since the Pontryagin classes are torsion-free, it is enough to verify the thesis of the lemma in the de'Rham cohomologies. Let us choose a basis $(E_1,E_2,E_3)$ of the Lie algebra $\soa(3)$ satisfying the standard commutation relations $[E_1,E_2]=E_3,\ldots cycl$. The representations $\rho_3,\rho_5$ maps this basis to the following matrices (see \cite{bn}):
\begin{gather*}
  \rho_3(E_1)=\left(\begin{smallmatrix}0&0&0\\ 0&0&1\\ 0&\scm 1&0\\ \end{smallmatrix}\right), \quad \rho_3(E_2)=\left(\begin{smallmatrix}0&0&1\\ 0&0&0\\ \scm 1&0&0\\ \end{smallmatrix}\right), \quad \rho_3(E_3)=\left(\begin{smallmatrix}0&1&0\\ \scm 1&0&0\\ 0&0&0\\ \end{smallmatrix}\right),\\
\rho_5(E_1)=\left(\begin{smallmatrix} 
0&0&0&0&\scriptscriptstyle{\sqrt{3}}\\
0&0&1&0&0\\
0&\scm1&0&0&0\\
0&0&0&0&1\\
\scriptscriptstyle{-\sqrt{3}}&0&0&\scm1&0
\end{smallmatrix}\right),\quad
\rho_5(E_2)=\left(\begin{smallmatrix} 
0&0&\scriptscriptstyle{\sqrt{3}}&0&0\\
0&0&0&0&1\\
\scriptscriptstyle{-\sqrt{3}}&0&0&1&0\\
0&0&\scm1&0&0\\
0&\scm1&0&0&0
\end{smallmatrix}\right),\quad
\rho_5(E_3)=\left(\begin{smallmatrix} 
0&0&0&0&0\\
0&0&0&2&0\\
0&0&0&0&1\\
0&\scm2&0&0&0\\
0&0&\scm1&0&0
\end{smallmatrix}\right).
\end{gather*}
The following identities hold
\begin{gather}
  \label{chpol}
%\begin{split}
  \det\big(\lambda\; \mathrm{I} +r_1 \rho_3(E_1)+ r_2 \rho_3(E_2) +r_3 \rho_3(E_3)\big) = \lambda^3 + \lambda\; (r_1^2+r_2^2+r_3^2),\\
  \det\big(\lambda\; \mathrm{I} +r_1 \rho_5(E_1)+ r_2 \rho_5(E_2) +r_3 \rho_5(E_3)\big) = \lambda^5 + \lambda^3\; 5\cdot(r_1^2+r_2^2+r_3^2) + \lambda\; (\ldots).\nonumber
%4(r_1^2+r_2^2+r_3^2)^2.
\end{gather}

We choose an $\soa(3)$ connection $\Gamma$ on the principal bundle $\pb$. The local curvature form of $\Gamma$ reads
\[
K = r^1 E_1+ r^2 E_2 +r^3 E_3,\qquad  r^j \in \Omega^2.
\]
The differential form representing the Pontrjagin class $p_1$ is constructed from the invariant polynomial of the curvature (see \cite{milch}). Using the identities (\ref{chpol}) we get the relation (\ref{ponrel}).\\
\end{proof}

\end{document}